\documentclass[11pt]{amsart}

\usepackage{amsthm}
\usepackage{amssymb}
\usepackage{amsmath}

\newcommand{\labell}[1] {\label{#1}}

\numberwithin{equation}{section}

\newtheorem {Theorem}   {Theorem} 
\numberwithin{Theorem}{section}

\newtheorem {Proposition}[Theorem]{Proposition}  
\theoremstyle{definition}

\theoremstyle{remark}
\newtheorem{Remark}[Theorem]{Remark}
\newtheorem{Example}[Theorem]{Example}
\newtheorem{Conjecture}[Theorem]{Conjecture}

\newtheorem {Corollary}[Theorem]{Corollary}

\newcommand{\CH}{{\mathcal H}}

\newcommand{\AP}{{\mathcal A}{\mathcal P}}

\def	\C {{\mathbb C}}
\def	\reals	{{\mathbb R}}
\def	\R	{{\mathbb R}}

\def	\TT	{{\mathbb T}}

\def	\CP	{{\mathbb C}{\mathbb P}}

\def	\irr	{\operatorname{irr}}
\def	\HZ	{\scriptscriptstyle\operatorname{HZ}}
\def	\SB	{\operatorname{SB}}
\def	\CL	{\operatorname{CL}}

\begin{document}




\title[The Hamiltonian Seifert Conjecture]{
The Hamiltonian Seifert Conjecture: Examples and Open Problems}

\author[Viktor Ginzburg]{Viktor L. Ginzburg}

\address{Department of Mathematics, 
         UC Santa Cruz, 
         Santa Cruz, CA 95064, USA}
\email{ginzburg@math.ucsc.edu}

\date{March, 2000. This paper is submitted to the Proceedings of the
Third ECM}

\thanks{The work is partially supported by the NSF and by the faculty
research funds of the University of California, Santa Cruz.}


\begin{abstract}
Hamiltonian dynamical systems tend to have infinitely many
periodic orbits. For example, for a broad class of symplectic
manifolds almost all levels of a proper smooth Hamiltonian carry
periodic orbits. The Hamiltonian Seifert conjecture is the
existence problem for regular compact energy levels without periodic 
orbits.

Very little is known about how large the set of regular energy values 
without periodic orbits can be. For instance, in all known examples
of Hamiltonian flows on  linear spaces such energy values 
form a discrete set, whereas ``almost existence theorems'' only require 
this set to have zero measure. We describe constructions of
Hamiltonian flows without periodic orbits on one energy level and formulate 
conjectures and open problems.
\end{abstract}

\maketitle

\section{Introduction} 
Hamiltonian flows of proper smooth functions on $\R^{2n}$ or on many
other symplectic manifolds are known to have periodic orbits on almost all 
energy levels; see \cite{FHV,hv-hs,ho-ze:book}. On the other hand, there
exist examples of proper smooth Hamiltonians on $\R^{2n}$ with one regular 
level carrying no periodic orbits, \cite{gi:seifert, gi:seifert97, 
herman-fax}. 

In this notice we address the question of how large the set
$\AP_H$ of regular values of $H$ without periodic orbits can be.
We refer to this problem as the Hamiltonian Seifert conjecture
(see \cite{kuper-k-cong,kuper-k}). Below
we recall some results and notions
concerning existence of periodic orbits and constructions of examples
without periodic orbits. We also formulate conjectures
on the size of the aperiodic set $\AP_H$. Very little is known on this 
problem and the conjectures are made exclusively based on 
supporting
evidence rather than on the lack of examples. 

Another subject we
briefly touch upon is Hamiltonian systems describing the motion of
a charge in a magnetic field. (See \cite{gi:Cambr} for an introduction
and a survey of results.) These systems play an important role in the
constructions of ``counterexamples'' to the Hamiltonian Seifert 
conjecture.
A detailed review of these counterexamples in the
context of dynamical systems can be found in \cite{gi:survey99}.

Two relevant questions are entirely omitted from our discussion: 
the existence of periodic orbits on contact type energy levels
(see, e.g., \cite{FHV,hv-hs,ho-ze:book,LT,vi:Theorem})
and the dynamics of the transition from low to high energy levels
for twisted geodesic flows (see, e.g., \cite{PP} and \cite{gi:Cambr}).

\bigskip\noindent\textbf{Acknowledgments.} The author is grateful to
Yael Karshon and Ely Kerman for their advice and numerous useful 
suggestions.

\section[Existence Results]{General Existence Results for Periodic 
Orbits of Hamiltonian Systems}
\labell{sec:existence}

Let $W$ be a symplectic manifold and $H\colon W\to \R$ a proper smooth 
Hamiltonian. One cannot 
guarantee that a given regular level of $H$ carries a periodic orbit
even when $W=\R^{2n}$, as we will see in Section \ref{sec:seifert}.
However, for a broad class of symplectic manifolds it is true 
that periodic orbits exist on almost all levels of $H$. (Here
``almost all'' is understood in the sense of measure theory.)

\begin{Theorem}[Almost Existence of Periodic Orbits]
\labell{thm:almost-exist}
Let $P$ be a compact symplectic manifold with $\pi_2(P)=0$. Then
for $W=P\times \R^{2l}$ with its split symplectic structure, almost
all energy levels of a proper smooth Hamiltonian
$H\colon W\to \R$ carry periodic orbits.
\end{Theorem}

This theorem is proved for $W=\R^{2n}$ by Hofer, Zehnder, and Struwe
(see \cite{ho-ze:per-sol,ho-ze:book,str}) and in the general case
by Floer, Hofer, and 
Viterbo, \cite{FHV}.

Analogues of Theorem \ref{thm:almost-exist} hold for many other
manifolds $W$, for example, for $W=\CP^n$, \cite{hv-hs}. 
(See also \cite{LT}.) The theorem also 
holds for $W=P\times D^2(r)$, where $\pi_2(P)$ is not necessarily trivial,
provided that $H$ is proper as a function to its image and that the 
radius $r>0$ is less than a certain constant $m(P,\omega)$ 
depending on the symplectic topology of $P$, \cite{hv-hs}.

The Hofer--Zehnder capacity $c_{\HZ}$ provides
a convenient setting for working with the almost existence problem.
Let $(W,\omega)$ be a symplectic manifold. Consider the class $\CH$ 
of smooth non-negative functions on $W$ such that

\begin{enumerate}

\item a function $H\in \CH$ vanishes on some open set, depending on $H$;

\item a function $H\in \CH$ assumes its maximal value $m(H)$ on a 
complement to a compact set, which again may depend on $H$; and

\item the Hamiltonian flow of $H\in \CH$ does not have non-constant
periodic orbits of period $T\leq 1$. 

\end{enumerate}
By definition, $c_{\HZ}(W,\omega)=\sup m(H)$ over all $H\in \CH$. 
Clearly, $c_{\HZ}$ takes values in $(0,\infty]$. The
capacity $c_{\HZ}$ is monotone with respect to inclusions and
homogeneous of  degree one with respect to conformal symplectic 
diffeomorphisms 
(in contrast with the symplectic volume, which is homogeneous of 
degree $n$, where $2n$ is the dimension). In addition, it satisfies
the following normalization condition:
$
c_{\HZ}(D^{2n}(r))=c_{\HZ}(D^{2}(r)\times \R^{2n-2})=\pi r^2,
$
where $D^{2k}(r)$ is the ball of radius $r$ in $\R^{2k}$ equipped with
the standard symplectic structure. (The proof of the normalization
condition is non-trivial. See \cite{ho-ze:book} for more details and
further references.) 

In what follows, we say that $W$ has bounded capacity if
for every open set $U\subset W$ with compact closure, 
$c_{\HZ}(U)<\infty$. Note that a non-compact manifold $W$
with $c_{\HZ}(W)=\infty$ can have bounded capacity (e.g., $W=\R^{2n}$). 
However, 
$c_{\HZ}(W)<\infty$ does imply, by monotonicity, that $W$ has bounded 
capacity.

\begin{Proposition}[Hofer--Zehnder, \cite{ho-ze:book}]
\labell{prop:capacity}
Assume that $W$ has bounded capacity. Then the almost existence
theorem holds for every smooth function $H$ on $W$ such that
the sets $\{H\leq a\}$ are compact.
\end{Proposition}

This proposition is a relatively straightforward consequence of the
definition  of the Hofer--Zehnder capacity, 
see \cite[Section 4.2]{ho-ze:book}. However, proving that $W$ has 
bounded capacity essentially amounts to showing that periodic orbits exist 
on a dense set of levels, which is quite non-trivial already for
$\R^{2n}$. The manifolds
$W=P\times \R^{2n}$ (with $\pi_2(P)=0$), $\CP^n$,
and also $P\times D^2(r)$ (with
$r<m(P,\omega)$) have bounded capacity, \cite{FHV,ho-ze:per-sol,
ho-ze:book,hv-hs}), as do all manifolds for which the almost existence 
theorem has been
established. (It seems to be unknown whether or not $\TT^{2n}$ with the 
standard symplectic structure has finite capacity for
$n>1$. However, almost existence 
has been proven for $H\colon \TT^{2n}\to \R$ satisfying some additional 
conditions; see, e.g., \cite{jiang}.)

Almost existence of periodic orbits
(Theorem \ref{thm:almost-exist}) does not extend to all compact 
symplectic manifolds:

\begin{Example}[Zehnder's torus, \cite{zehnder}]
\labell{ex:Zehn}
Let $2n\geq 4$. Consider the torus $W=\TT^{2n}$ with an
irrational translation-invariant symplectic structure $\omega_{\irr}$.
Choose a Hamiltonian $H$ on $W$ so that every level $\{ H=c\}$ with
$c\in (0.5, 1.5)$ is the union of two standard embedded tori 
$\TT^{2n-1}\subset \TT^{2n}$. Since $\omega_{\irr}$
is irrational, the Hamiltonian flow of $H$ on $\TT^{2n-1}$ is
quasiperiodic. Thus, none of the levels
$\{H=c\}$ with $c\in (0.5, 1.5)$ carries a periodic orbit. 

According to M. Herman, \cite{herman1,herman2}, this phenomenon is stable:
the flow of a sufficiently $C^{2n+\delta}$-small perturbation of $H$ 
still has no periodic orbits in the energy range $(0.6, 1.4)$, provided
that $\omega_{\irr}$ satisfies a certain Diophantine condition.
\end{Example}

As a consequence of Proposition \ref{prop:capacity}, 
$c_{\HZ}(\TT^{2n},\omega_{\irr})=\infty$, even though $\TT^{2n}$ is 
compact. Note also that in
Zehnder's example the symplectic form is not exact on the energy levels
without periodic orbits. However, many flows for which almost existence
has been established (e.g., on $P\times \R^2$) also have this property.
Thus, although it feels that non-existence of periodic orbits
in Zehnder's example is forced by the topology of $\omega_{\irr}$, it is
not clear how to make rigorous sense of this statement.

Now we are in a position to state the main question considered in
this paper. We call a value $a$ of $H$ \emph{aperiodic} if the level
$\{H=a\}$ carries no periodic orbits.
\begin{quote}
Let $W$ be a symplectic manifold of bounded capacity and $H$ a smooth
proper function on $W$. How large can the set $\AP_H$ of regular 
aperiodic values of $H$ be?
\end{quote}
As we show in the next section, the set $\AP_H$ can be non-empty 
for an appropriately chosen function $H$ on $\R^{2n}$, and hence,
by Darboux's theorem, on any symplectic manifold.

\section{The Hamiltonian Seifert Conjecture} 
\labell{sec:seifert}

\subsection{Results}
The \emph{Seifert conjecture} is the question, posed by Seifert in 
\cite{Seifert-1950}, whether or not every non-vanishing vector field on 
$S^3$ has a periodic orbit. Of course, in this question the sphere $S^3$
can be replaced by another manifold, and additional restrictions can
be imposed on the vector field. For example, the Seifert conjecture 
can be formulated for $C^1$- or $C^\infty$-smooth or real analytic
vector fields, divergence--free vector fields, etc. Note that the 
Seifert conjecture, when interpreted as above, is a question 
rather than a conjecture. Yet, we will refer to negative answers to 
this question as counterexamples to the Seifert conjecture.

The Seifert conjecture has a rich history extending for over than forty 
years. Here we mention only some of the relevant results. The first 
breakthrough was due to Wilson, \cite{wilson}, who 
constructed a $C^\infty$-smooth vector field without periodic orbits
on $S^{2n+1}$, $2n+1\geq 5$. A $C^1$-smooth non-vanishing
vector field without periodic
orbits on $S^3$ was found by Schweitzer, \cite{schweitzer}, and a 
$C^2$-vector field by Harrison, \cite{harrison}. Finally, 
a real analytic non-vanishing vector field on $S^3$ without periodic 
orbits was constructed by K. Kuperberg, \cite{kuk}. The reader 
interested in a detailed
discussion and further references should consult 
\cite{gi:survey99,kug,kugk,kuper-k-cong, kuper-k}.

For Hamiltonian flows, the Seifert conjecture can be formulated in
a number of ways; see \cite{gi:survey99}. For example, one may ask if 
there is a proper smooth function on a given symplectic manifold 
(e.g., $\reals^{2n})$, having a regular aperiodic value,
or, more generally, as in the previous section, if there exists a
smooth proper function $H$ with a given set $\AP_H$. Here we state and 
discuss known results and make some conjectures.

Recall that a \emph{characteristic} of a two-form $\eta$ of 
rank $(2n-2)$ on a $(2n-1)$-dimensional manifold is an integral curve
of the line field formed by the null-spaces $\ker\eta$.
The Hamiltonian Seifert conjecture (in a form slightly weaker than
considered above) can be stated as
whether or not in a given symplectic manifold there exists
a regular compact hypersurface without closed characteristics. 

Let $i\colon M\hookrightarrow W$ be an embedded smooth compact 
hypersurface without boundary in a $2n$-dimensional symplectic manifold 
$(W,\omega)$. 

\begin{Theorem}[\cite{gi:seifert, gi:seifert97, herman-fax}]
\labell{Theorem:main}
Assume that $2n\geq 6$ and that $i^*\omega$ has a finite number of 
closed characteristics. Then there exists a $C^\infty$-embedding 
$i'\colon M\hookrightarrow W$, which is $C^0$-close and isotopic to $i$, 
such that ${i'}^*\omega$ has no closed characteristics.
\end{Theorem}

An irrational ellipsoid $M$ in the standard symplectic vector space 
$\C^{n}$ is the unit energy level of the quadratic
Hamiltonian $H_0=\sum \lambda_j|z_j|^2$, where the eigenvalues 
$\lambda_j$ are positive and rationally independent. 
The level $M$ carries exactly $n$ periodic orbits. 
Applying Theorem \ref{Theorem:main}, we obtain the following

\begin{Corollary}
\labell{Corollary:function}
For $2n\geq 6$, there exists a $C^\infty$-function 
$H\colon\reals^{2n}\to \reals$, $C^0$-close and isotopic
(with a compact support) to $H_0$,
such that the Hamiltonian flow of $H$ has no closed trajectories on 
the level set $\{ H=1\}$.
\end{Corollary}

\begin{Remark}
The method used in \cite{gi:seifert, gi:seifert97}
to construct such a function $H$ relies on a Hamiltonian version of 
Wilson's plug. The horocycle flow (see Example \ref{example:horocycle})
is employed to interrupt periodic orbits of the flow of $H_0$ on the
level $\{H_0=1\}$. The function $H$ differs from $H_0$ only in $n$ small
balls with centers on these periodic orbits. 
A simple outline of the proof of Theorem \ref{Theorem:main}, some
generalizations, and the list of all known Hamiltonian flows with
$\AP_H\neq\emptyset$
can be found in \cite{gi:survey99}. Here we just mention that
irrational ellipsoids are not the only hypersurfaces in $\R^{2n}$ with
a finite number of closed characteristics. Examples of such non-simply
connected hypersurfaces are found by Laudenbach, \cite{laud}.

The construction given in \cite{herman-fax}  uses plugs arising
as symplectizations of Wilson's plug for $2n\geq 8$ and Harrison's  
plug \cite{harrison} for $2n=6$. (Hence, the $C^{3-\epsilon}$-smoothness
of the flow obtained by this method in dimension six.)
\end{Remark}

\subsection{Conjectures and Open Problems}

Returning to the question of the possible size of the set $\AP_H$,
we primarily focus in this section on 
proper Hamiltonians $H$ on $\R^{2n}$. (See the next section and
\cite{gi:survey99} for results and conjectures for some other
symplectic manifolds.) Recall that an upper bound on 
this set is given by Theorem \ref{thm:almost-exist}: the set 
$\AP_H$ 
must have zero measure. By Theorem \ref{Theorem:main} and Corollary
\ref{Corollary:function}, $\AP_H$ can be non-empty.

There is almost no doubt that the method used in 
\cite{gi:seifert, gi:seifert97} to modify the function $H_0$
can be applied, literally without any change, to a sequence of energy
levels $a_k$ converging to some value $a_0$. The resulting function
$H$ can then be smoothened along the level $\{H=a_0\}$ at the 
expense of turning $a_0$ into a critical value of infinite order 
degeneration. Thus, $\AP_H$ can contain a sequence 
converging to a singular, infinitely degenerate, value of $H$.

The next step is more subtle. It appears plausible that the function
$H_0$ can be modified as in the proof of Theorem \ref{Theorem:main}
simultaneously for all energy values from a compact 
zero measure set $K\subset(0,\infty)$. The resulting function $H$ will 
then be smoothly isotopic to $H_0$ and have $K$ as the set of aperiodic 
values. This leads to the following

\begin{Conjecture}
For $2n\geq 6$, there exists a $C^\infty$-function 
$H\colon\reals^{2n}\to \reals$, $C^0$-close and isotopic
(with a compact support) to $H_0$,
such that the Hamiltonian flow of $H$ has no closed trajectories for
energy values in the set $K$.
\end{Conjecture}

\begin{Remark}
This method is unlikely to be applicable to the construction of a 
function $H$ with a dense set $\AP_H$. The difficulty can be best seen
by considering the function $H$ from Corollary \ref{Corollary:function}.
The characteristic foliation of $\omega$ on the level $\{H=a\}$, 
where $a\neq 1$ and $a>0$, is diffeomorphic to that on $\{H_0=a\}$. 
(In fact, the forms $\omega|_{\{H=a\}}$ and $\omega|_{\{H_0=a\}}$ 
are conformally equivalent.) As a consequence,
the flow of $H$ on $\omega|_{\{H=a\}}$ has exactly $n$ periodic orbits
and these orbits are non-degenerate and hence stable under small
perturbations. The ``bifurcation'' that happens at $a=1$ is simply 
that the periods of these orbits go to infinity as $a\to 1$. Other 
systems obtained by this method from any function $H_0$ with 
non-degenerate orbits will exhibit a similar behavior: 
periodic orbits on $\{H=a\}$, where $a\not\in \AP_H$, 
will be non-degenerate. On the other hand, for a system with a dense
set $\AP_H$ every energy level must have only degenerate periodic orbits.
\end{Remark}

\begin{Remark}
Neither the method from \cite{gi:seifert, gi:seifert97} nor
from \cite{herman-fax} apply to four-dimensional manifolds. There are also
serious difficulties in adapting Kuperberg's plug, \cite{kuk},
to the symplectic setting. 
A $C^1$-smooth divergence--free vector field on $S^3$ was
found by G. Kuperberg, \cite{kug}. It is not known whether or not
this vector field can be obtained by a $C^2$- (or even $C^1$-) embedding
of $S^3$ into $\R^4$. Such an 
embedding would give a $C^2$-smooth ``counterexample'' to the Hamiltonian
Seifert conjecture in dimension four.
\end{Remark}

\section{A Charge in a Magnetic Field}
\labell{sec:magnetic}
The first example of a Hamiltonian flow with a compact energy level
carrying no periodic orbits was found by Hedlund in 1936 when he
proved that the horocycle flow is minimal, i.e., every integral curve is
dense, \cite{hedlund}. To put Hedlund's result in the context of 
the Hamiltonian Seifert conjecture, consider a compact Riemannian manifold
$M$ equipped with a closed two-form $\sigma$. Denote by $\omega_0$ the
standard symplectic form on $T^*M$ and by $\pi\colon T^*M\to M$ the 
natural projection. The form $\omega=\omega_0+\pi^*\sigma$ 
on $T^*M$ is symplectic (a twisted symplectic structure). 
The flow of the standard metric Hamiltonian $H$ (the kinetic energy) 
on $(T^*M,\omega)$ describes the motion of a charge on
$M$ in the magnetic field $\sigma$. We will refer to this flow as a
\emph{twisted geodesic flow}. The behavior of twisted geodesic flows 
seems to depend on whether $\sigma$ is allowed to degenerate or 
not. In what follows we focus entirely on the case where $\sigma$ is 
non-degenerate, i.e., symplectic.

The study of twisted geodesic flows by methods of symplectic topology
originates from \cite{arnold}, where V. Arnold showed that for $M=\TT^2$
with a flat metric and a non-vanishing form $\sigma$ every level of
$H$ carries periodic orbits. This theorem is extended to low energy levels
on other compact surfaces in \cite{gi:FA,gi:MathZ}. (See \cite{gi:Cambr}
for a survey of related results.) The existence of periodic orbits on
low energy levels has also been established for many manifolds $M$,
when $2m=\dim M>2$. For example, when $\sigma$ is symplectic and 
compatible with the metric on $M$, every low energy level carries at least
$\CL(M)+m$ periodic orbits, \cite{ely:paper}, and at least $\SB(M)$ 
periodic orbits when the orbits are non-degenerate, \cite{gi-ke}. 
(Here $\CL(M)$ is the cup--length of $M$ and $\SB(M)$ is the sum of 
Betti numbers of $M$.) Note also that this problem of existence of
periodic orbits can be extended to a broader class of Hamiltonian
flows to include the Weinstein--Moser theorem, \cite{gi-ke,ely:paper}.

The following example shows that the above results do not generalize
to all energy levels.

\begin{Example}[The horocycle flow]
\labell{example:horocycle}
Let $M$ be a compact surface of genus $g\geq 2$ equipped with a
metric of constant negative curvature $K=-1$. Let $\sigma$ be the
area form on $M$. Then the Hamiltonian flow on the energy level
$\{H=1\}$ has no periodic orbits. In fact, the flow on this energy
level is smoothly equivalent to the horocycle flow (see, e.g., 
\cite{gi:MathZ}), which is minimal by a theorem of Hedlund, 
\cite{hedlund}. This example has (Hamiltonian) analogues in
all even dimensions; see \cite[Example 4.2]{gi:survey99}.
\end{Example}

It is also known that a neighborhood of $M$ in $T^*M$ with the
twisted symplectic structure has finite Hofer--Zehnder capacity 
under some additional assumptions on $M$, \cite{Lu}. Moreover, 
for some manifolds $M$, the twisted cotangent bundle 
has bounded capacity, \cite{gi-ke,jiang,Lu}. This implies almost
existence of periodic orbits for a certain class of twisted geodesic 
flows (e.g., on $\TT^{2m}$ with any form $\sigma$). 

\begin{Conjecture}
\labell{conj:magn}
For every $M$ and any symplectic form $\sigma$, a neighborhood
of $M$ in $T^*M$ has finite Hofer--Zehnder capacity. Moreover, 
almost all levels of $H$ carry periodic orbits. When $\sigma$ is
symplectic, every low energy level of $H$ has a contractible periodic 
orbit.

\end{Conjecture}

\begin{Remark}
This conjecture is supported by some additional evidence; see, e.g.,
\cite{po}. Interestingly, there seems to be no sufficient evidence 
that $T^*M$ has bounded Hofer--Zehnder capacity in general. For example, 
it is 
not known whether in the setting of Example \ref{example:horocycle} 
the sets $\{H<a\}$ with $a>1$ have finite capacity. 
\end{Remark}

Example \ref{example:horocycle}  is the only known ``naturally arising''
example of a Hamiltonian flow with a regular compact aperiodic energy 
level. Furthermore, this is the only known example of a twisted geodesic 
flow with an aperiodic energy level. It is not known if there exists 
a twisted geodesic flow with more than one aperiodic energy level.

\end{document}